\documentclass[a4paper,12pt]{article}
\usepackage{amsfonts}
\usepackage{mathrsfs}
\usepackage{algorithm, algpseudocode}
\usepackage{amsmath,amssymb,amsthm,latexsym,amsfonts}
\usepackage{pstricks}
\usepackage{enumerate}

\usepackage{graphicx}
\usepackage{color}
\usepackage{ifpdf}
\usepackage{epstopdf}
\usepackage{caption}

\begin{document}\baselineskip 0.8cm

\newtheorem{lem}{Lemma}[section]
\newtheorem{thm}[lem]{Theorem}
\newtheorem{cor}[lem]{Corollary}
\newtheorem{exa}[lem]{Example}
\newtheorem{con}[lem]{Conjecture}
\newtheorem{rem}[lem]{Remark}
\newtheorem{obs}[lem]{Observation}
\newtheorem{definition}[lem]{Definition}
\newtheorem{pro}[lem]{Proposition}
\newtheorem{prob}[lem]{Problem}
\theoremstyle{plain}
\newcommand{\D}{\displaystyle}
\newcommand{\DF}[2]{\D\frac{#1}{#2}}

\renewcommand{\figurename}{{\bf Fig}}
\captionsetup{labelfont=bf}

\title{A Survey on Monochromatic Connections of Graphs\footnote{Supported by NSFC No.11371205 and 11531011.}}
\author{\small Xueliang Li, Di Wu\\
      {\it\small Center for Combinatorics and LPMC}\\
      {\it\small Nankai University, Tianjin 300071, P.R. China}\\
       {\it\small lxl@nankai.edu.cn; wudiol@mail.nankai.edu.cn}}
\date{}
\maketitle

\begin{abstract}
The concept of monochromatic connection of graphs was introduced by Caro and Yuster in 2011.
Recently, a lot of results have been published about it.
In this survey, we attempt to bring together all the results that dealt with it.
We begin with an introduction, and then classify the results into the following categories:
monochromatic connection coloring of edge-version, monochromatic connection coloring of vertex-version,
monochromatic index, monochromatic connection coloring of total-version.
\\[2mm]


\end{abstract}

\section{Introduction}

All graphs considered in this paper are simple, finite and
undirected. We follow the terminology and notation of Bondy and
Murty in \cite{Bondy}.
Let $H$ be a nontrivial connected graph with an edge-coloring $f:
E(H)\rightarrow \{1,2,\ldots,\ell\}$ ($\ell$ is a positive integer),
where adjacent edges may be colored the same.
A path in an edge-colored graph $H$ is a \emph{monochromatic path}
if all the edges of the path are colored with a same color.
The graph $H$ is called \emph{monochromatically connected},
if any two vertices of $H$ are connected by a monochromatic path.
An edge-coloring of $H$ is a \emph{monochromatic connection coloring}
(\emph{MC-coloring}) if it makes $H$ monochromatically connected.
How colorful can an MC-coloring be?
This question is the natural opposite of the well-studied problem of
rainbow connection coloring \cite{Caro1, Chartrand1, Krivelevich and Yuster, Sun, LiSun2},
where in the latter we seek to find an edge-coloring with minimum number of colors
so that there is a rainbow path joining any two vertices; see \cite{Sun, LiSun2} for details.
As introduced by Caro and Yuster \cite{Caro},
for a connected graph $G$, the \emph{monochromatic connection number} of $G$, denoted by $mc(G)$,
is the maximum number of colors that are needed in order to make $G$
monochromatically connected.
An \emph{extremal MC-coloring} is an MC-coloring that uses $mc(G)$ colors.

Gonz¨¢lez-Moreno, Guevara, and Montellano-Ballesteros\cite{Gonz¨¢lez} generalized the above concept to digraphs.
Let $H$ be a nontrivial strongly connected digraph with an arc-coloring $f:
A(H)\rightarrow \{1,2,\ldots,\ell\}$ ($\ell$ is a positive integer), where adjacent arcs may be colored the same.
A path in an arc-colored graph $H$ is a \emph{monochromatic path}
if all the arcs of the path are colored with a same color.
The strongly connected digraph $H$ is called \emph{strongly monochromatically connected}£¬
if for every pair $\{u,v\}$ of vertices in $H$,
there exist both a $(u,v)$-monochromatic path and a $(v,u)$-monochromatic path.
An arc-coloring of $H$ is a \emph{strongly monochromatic connection coloring}
(\emph{SMC-coloring}) if it makes $H$ strongly monochromatically connected.
For a strongly connected digraph $D$, the \emph{strongly monochromatic connection number} of $D$, denoted by $smc(D)$,
is the maximum number of colors that are needed in order to make $D$
strongly monochromatically connected.
An \emph{extremal SMC-coloring} is an SMC-coloring that uses $smc(D)$ colors.

Now we introduce another generalization of the monochromatic connection number by Li and Wu \cite{Di1}.
A tree $T$ in an edge-colored graph $H$ is called a \emph{monochromatic tree} if all the edges of $T$ have the same color. For an $S\subseteq V(H)$,
a \emph{monochromatic $S$-tree} is a monochromatic tree of $H$ containing the vertices of $S$.
Given an integer $k$ with $2\leq k\leq |V(H)|$, the graph $H$ is called \emph{$k$-monochromatically connected}
if for any set $S$ of $k$ vertices of $H$, there exists a monochromatic $S$-tree in $H$.
An edge-coloring of $H$ is called a \emph{$k$-monochromatic connection coloring} (\emph{$MX_k$-coloring}) if it makes $H$ $k$-monochromatically connected.
For a connected graph $G$ and a given integer $k$ such that $2\leq k\leq |V(G)|$,
the \emph{$k$-monochromatic index $mx_k(G)$} of $G$ is the maximum number of
colors that are needed in order to make $G$ $k$-monochromatically connected.
An \emph{extremal $MX_k$-coloring} is an $MX_k$-coloring that uses $mx_k(G)$ colors.
By definition, we have $mx_{|V(G)|}(G)\leq \cdots\leq mx_3(G)\leq mx_2(G)=mc(G)$.

Note that the above graph-parameters are defined on edge-colored graphs.
Naturally,  Cai, Li and Wu \cite{Di2} introduced a graph-parameter
corresponding to monochromatic connection number which is defined on vertex-colored graphs.
Let $H$ be a nontrivial connected graph with a vertex-coloring $f:
V(H)\rightarrow \{1,2,\ldots,\ell\}$ ($\ell$ is a positive integer), where adjacent vertices may be colored the same.
A path in a vertex-colored graph $H$ is a \emph{vertex-monochromatic path}
if all the internal vertices of the path are colored with a same color.
The graph $H$ is called \emph{vertex-monochromatically connected},
if any two vertices of $H$ are connected by a vertex-monochromatic path.
A vertex-coloring of $H$ is a \emph{vertex-monochromatic connection coloring}
(\emph{VMC-coloring}) if it makes $H$ vertex-monochromatically connected.
For a connected graph $G$, the \emph{vertex-monochromatic connection number} of $G$, denoted by $vmc(G)$,
is the maximum number of colors that are needed in order to make $G$
vertex-monochromatically connected.
An \emph{extremal VMC-coloring} is a VMC-coloring that uses $vmc(G)$ colors.

Li and Wu \cite{Di1} introduced another graph-parameter corresponding to the $k$-mono- chromatic index, which is defined on vertex-colored graphs.
A tree $T$ in a vertex-colored graph $H$ is called a \emph{vertex-monochromatic tree}
if all the internal vertices of $T$ have the same color. For an $S\subseteq V(H)$,
a \emph{vertex-monochromatic $S$-tree} is a vertex-monochromatic tree of $H$ containing the vertices of $S$.
Given an integer $k$ with $2\leq k\leq |V(H)|$, the graph $H$ is called \emph{$k$-vertex-monochromatically connected}
if for any set $S$ of $k$ vertices of $H$, there exists a vertex-monochromatic $S$-tree in $H$.
A vertex-coloring of $H$ is called a \emph{$k$-vertex-monochromatic connection coloring} (\emph{$VMX_k$-coloring})
if it makes $H$ $k$-vertex-monochromatically connected.
For a connected graph $G$ and a given integer $k$ such that $2\leq k\leq |V(G)|$,
the \emph{$k$-vertex-monochromatic index $vmx_k(G)$} of $G$ is the maximum number of
colors that are needed in order to make $G$ $k$-vertex-monochromatically connected.
An \emph{extremal $VMX_k$-coloring} is a $VMX_k$-coloring that uses $vmx_k(G)$ colors.
By definition, we have $vmx_{|V(G)|}(G)\leq \cdots\leq vmx_3(G)\leq vmx_2(G)=vmc(G)$.

Jiang, Li and Zhang \cite{Jiang1} introduced the monochromatic connection of total-coloring version.
Let $H$ be a nontrivial connected graph with a total-coloring $f:
V(H)\cup E(H)\rightarrow \{1,2,\ldots,\ell\}$ ($\ell$ is a positive integer),
where any two elements may be colored the same.
A path in a total-colored graph $H$ is a \emph{total-monochromatic path}
if all the edges and internal vertices of the path are colored with a same color.
The graph $H$ is called \emph{total-monochromatically connected}£¬
if any two vertices of $H$ are connected by a total-monochromatic path.
A total-coloring of $H$ is a \emph{total-monochromatic connection coloring}
(\emph{TMC-coloring}) if it makes $H$ total-monochromatically connected.
For a connected graph $G$, the \emph{total-monochromatic connection number} of $G$, denoted by $tmc(G)$,
is the maximum number of colors that are needed in order to make $G$
total-monochromatically connected.
An \emph{extremal TMC-coloring} is a TMC-coloring that uses $tmc(G)$ colors.

Next, we recall the definitions of various products of graphs, which will be used in the sequel.
The \emph{Cartesian product} of two graphs $G$ and $H$, denoted
by $G\Box H$, is defined to have the vertex-set $V(G)\times V(H)$, in
which two vertices $(g,h)$ and $(g',h')$ are adjacent if and only if
$g=g'$ and $hh'\in E(H)$, or $h=h'$ and $gg'\in E(G)$.
The \emph{lexicographic product} $G\circ H$ of two graphs $G$ and $H$
has the vertex-set $V(G\circ H)=V(G)\times V(H)$, and two vertices
$(g,h),(g',h')$ are adjacent if and only if $gg'\in E(G)$, or $g=g'$ and
$hh'\in E(H)$.
The \emph{strong product} $G\boxtimes H$ of two graphs $G$ and $H$ has
the vertex-set $V(G)\times V(H)$. Two vertices $(g,h)$ and $(g',h')$
are adjacent whenever $gg'\in E(G)$ and $h=h'$, or $g=g'$ and $hh'
\in E(H)$, or $gg'\in E(G)$ and $hh'\in E(H)$.
The \emph{direct product} $G\times H$ of two graphs $G$ and $H$ has the
vertex-set $V(G)\times V(H)$. Two vertices $(g,h)$ and $(g',h')$ are
adjacent if the projections on both coordinates are adjacent, i.e.,
$gg'\in E(G)$ and $hh'\in E(H)$. Finally, the {\it join} $G+H$ of two graphs
$G$ and $H$ has the vertex-set $V(G)\cup V(H)$ and edge-set
$E(G)\cup E(H) \cup \{xy \ | \ x\in V(G), \ y\in V(H)\}$.

The most frequently occurring probability models of
random graphs is the Erd\"{o}s-R\'{e}nyi random graph model $G(n,p)$ \cite{ER}.
The model $G(n,p)$ consists of all graphs with
$n$ vertices in which the edges are chosen independently and with
probability $p$. We say an event $\mathcal{A}$ happens
\textit{with high  probability} if the probability that it happens approaches
$1$ as $n\rightarrow \infty $, i.e., $Pr[\mathcal{A}]=1-o_n(1)$. Sometimes,
we say \textit{w.h.p.} for short.
We will always assume that $n$ is the variable that tends to infinity.
Let $G$ and $H$ be two graphs on $n$ vertices. A property $P$ is said to be \emph{monotone} if
whenever $G\subseteq H$ and $G$ satisfies $P$, then $H$ also satisfies $P$. For a graph property $P$, a function $p(n)$ is called a threshold
function of $P$ if:
\begin{itemize}
\item for every $r(n) = \omega(p(n))$, $G(n, r(n))$ w.h.p. satisfies $P$; and

\item for every $r'(n) = o(p(n))$, $G(n, r'(n))$ w.h.p.
does not satisfy $P$.
\end{itemize}

Furthermore, $p(n)$ is called a sharp threshold function of $P$ if
there exist two positive constants $c$ and $C$ such that:
\begin{itemize}
  \item for every $r(n) \geq C\cdot p(n)$, $G(n, r(n))$ w.h.p. satisfies $P$; and
  \item for every $r'(n) \leq c\cdot p(n)$, $G(n, r'(n))$ w.h.p.
does not satisfy $P$.
\end{itemize}
It is well known that all monotone graph properties have a sharp
threshold function; see \cite{BT} and \cite{FK} for details.

\section{The edge-coloring version}

\subsection{Upper and lower bounds for $mc(G)$}

In \cite{Caro}, Caro and Yuster observed that a general lower bound for $mc(G)$ is $m(G)-n(G)+2$.
Simply color the edges of a spanning tree of $G$ with one color, and each of the
remaining edges with a distinct new color. Then, Caro and
Yuster gave some sufficient conditions for graphs attaining this
lower bound.
\begin{thm}\cite{Caro}
Let $G$ be a connected graph with $n>3$ vertices and $m$ edges.
If $G$ satisfies any of the following properties, then $mc(G)=m-n+2$.\\
$(a)$ $\overline{G}$ is 4-connected.\\
$(b)$ $G$ is triangle-free.\\
$(c)$ $\Delta(G)<n-\frac{2m-3(n-1)}{n-3}$.
In particular, this holds if $\Delta(G)\leq(n+1)/2$ or $\Delta(G)\leq n-2m/n$.\\
$(d)$ $diam(G)\geq3$.\\
$(e)$ $G$ has a cut vertex.
\end{thm}

Jin, Li and Wang got some conditions on graphs containing triangles.

\begin{thm}\cite{Jin}
Let $G$ be a connected graph of order $n \geq 7$. If $G$ does not have
subgraphs isomorphic to $K_4^-$, then $mc(G) = m -n + 2$, where $K_4^-$
denotes the graph obtained from $K_4$ by deleting an edge.
\end{thm}

\begin{thm}\cite{Jin}
Let $G$ be a connected graph of order $n \geq 7$. If $G$ does not have two
triangles that have exactly one common vertex, then $mc(G) = m - n + 2$.
\end{thm}

\begin{thm}\cite{Jin}
Let $G$ be a connected graph of order $n \geq 7$. If $G$ does not
have two vertex-disjoint triangles, then $mc(G) = m - n + 2$.
\end{thm}

Caro and Yuster \cite{Caro} also showed some nontrivial upper bounds for $mc(G)$ in terms of the chromatic number, the connectivity, and the
minimum degree. Recall that a graph is called \emph{$s$-perfectly-connected} if it can be partitioned into $s+1$
parts $\{v\},V_1,\ldots,V_s$, such that each $V_j$ induces a
connected subgraph, any pair $V_j,V_r$ induces a corresponding
complete bipartite graph, and $v$ has precisely one neighbor in each
$V_j$. Notice that such a graph has minimum degree $s$, and $v$ has
degree $s$.

\begin{thm}\cite{Caro}\label{Carothm1}

\noindent $(1)$ If $G$ is a complete $r$-partite graph, then $mc(G)=m-n+r$.\\
$(2)$ Any connected graph $G$ satisfies $mc(G)\leq m-n+\chi(G)$.\\
$(3)$ If $G$ is not $k$-connected, then $mc(G)\leq m-n+k$. This is sharp for any $k$.\\
$(4)$ If $\delta(G)=s$, then $mc(G)\leq m-n+s$, unless $G$ is $s$-perfectly-connected, in which case $mc(G)=m-n+s+1$.
\end{thm}

As an application of Theorem\ref{Carothm1}(4), Caro and Yuster got the upper bounds for the following planar graphs.
\begin{cor}\cite{Caro}

\noindent $(1)$ For $n\geq 5$, the wheel $G=W_n$ has $mc(G)=m-n+3$.\\
$(2)$ If $G$ is an outerplanar graph, then $mc(G)=m-n+2$, except that $mc(K_1\vee P_{n-1})=m-n+3$.\\
$(3)$ If $G$ is a planar graph with minimum degree 3, then $mc(G)\leq m-n+3$, except that $mc(K_2\vee P_{n-2})=m-n+4$.
\end{cor}

\subsection{Erd\H{o}s-Gallai-type problems for $mc(G)$}

Cai, Li and Wu \cite{Cai} studied the following two kinds of Erd\H{o}s-Gallai-type problems for $mc(G)$.

\noindent\textbf{Problem A.}
Given two positive integers $n$ and $k$ with $1\leq k\leq {n \choose 2}$, compute the minimum integer $f(n,k)$ such that
for any graph $G$ of order $n$, if $|E(G)|\geq f(n,k)$ then
$mc(G)\geq k$.

\noindent\textbf{Problem B.} Given two positive integers $n$ and $k$
with $1\leq k\leq {n \choose 2}$, compute the maximum integer $g(n,k)$ such
that for any graph $G$ of order $n$, if $|E(G)|\leq g(n,k)$ then
$mc(G)\leq k$.

It is worth mentioning that the two parameters $f(n,k)$ and $g(n,k)$
are equivalent to another two parameters. Let
$t(n,k)=\min\{|E(G)|:|V(G)|=n, mc(G)\geq k\}$ and
$s(n,k)=\max\{|E(G)|:|V(G)|=n, mc(G)\leq k\}$. It is easy to see that $t(n,k)=g(n,k-1)+1$ and $s(n,k)=f(n,k+1)-1$. In \cite{Cai}
the authors determined the exact values of $f(n,k)$ and $g(n,k)$ for all integers
$n$, $k$ with $1\leq k\leq {n \choose 2}$.

\begin{thm}\cite{Cai}
Given two positive integers $n$ and $k$ with $1\leq k\leq {n \choose 2}$,
\begin{displaymath}
f(n,k) = \left\{ \begin{array}{ll}
n+k-2 & \textrm{\ \ \ \ if\ $1\leq k\leq {n \choose 2}-2n+4$ \ \ }\\
{n \choose 2}+\left\lceil\frac{k-{n \choose 2}}{2}\right\rceil & \textrm{\ \ \ \
if\ ${n \choose 2}-2n+5\leq k\leq {n \choose 2}$.\ }
\end{array} \right.
\end{displaymath}
\end{thm}

\begin{thm}\cite{Cai}
Given two positive integers $n$ and $k$ with $1\leq k\leq {n \choose 2}$,
\begin{align*}
g(n,k)=
\begin{cases}
\binom{n}{2} & \textrm{\ \ \ \ if $k=\binom{n}{2}$\ }\\
k+t-1 & \textrm{\ \ \ \ if $\binom{n-t}{2}+t(n-t-1)+1\leq k\leq \binom{n-t}{2}+t(n-t)-1$ }\\
k+t-2 & \textrm{\ \ \ \ if $k=\binom{n-t}{2}+t(n-t)$}
\end{cases}
\end{align*}
for $2\leq t\leq n-1$.
\end{thm}

\subsection{Results for graph classes}

Gu, Li, Qin and Zhao\cite{Qin} characterized all connected graphs $G$ of size $m$ with small and large values of $mc(G)$.

\begin{thm}\cite{Qin}
Let $G$ be a connected graph. Then $mc(G)=1$ if and only if $G$ is a tree.
\end{thm}
\begin{thm}\cite{Qin}
Let $G$ be a connected graph. Then $mc(G)=2$ if and only if $G$ is a unicyclic graph except for $K_3$.
\end{thm}
\begin{thm}\cite{Qin}
Let $G$ be a connected graph. Then $mc(G)=3$ if and only if $G$ is either $K_3$ or a bicyclic graph except for $K_4-e$ .
\end{thm}
\begin{figure}[htbp]
\begin{center}
\includegraphics[scale = 0.7]{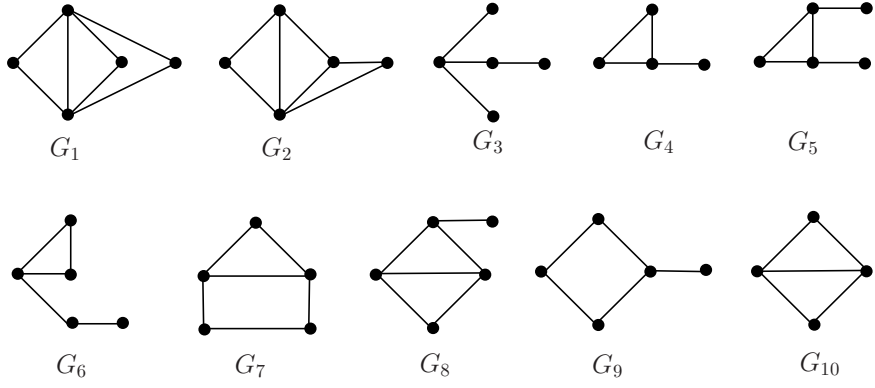}
\caption{The graphs in Theorem \ref{Qin-thm1} and Theorem \ref{Qin-thm2}}
\label{fig1}
\end{center}
\end{figure}

\begin{thm}\label{Qin-thm1}\cite{Qin}
Let $G$ be a connected graph. Then $mc(G)=4$ if and only if $G$ is either $K_4-e$ or a tricyclic graph except for $G_1,\ G_2,\ K_4$,
where $G_1,\ G_2$ are shown in Figure \emph{1}.
\end{thm}

\begin{thm}\cite{Qin}
Let $G$ be a connected graph. Then $mc(G)=m-1$ if and only if $G=K_n-e$.
\end{thm}

\begin{thm}\cite{Qin}
Let $G$ be a connected graph. Then $mc(G)=m-2$ if and only if $G\in \{K_n-2K_2,\ K_n-P_3,\ K_n-K_3,\ K_n-P_4\}$.
\end{thm}

\begin{thm}\cite{Qin}\label{Qin-thm2}
Let $G$ be a connected graph. Then $mc(G)=m-3$ if and only if $G\in \{K_n-3K_2,\ K_n-C_4,\ K_n-C_5,\ K_n-(P_2 \cup P_3),\ K_n-(P_2 \cup K_3), \ K_n-P_5,\ K_n-(P_2 \cup P_4), \ K_n-S_4,\ K_n-K_4,\ K_n-G_i\}$, where $i$ is an integer with $3 \leq i \leq 10$ and $G_i$ is shown in Figure \emph{1}.
\end{thm}

From the above theorems, they also verified the following corollary.
\begin{cor}\cite{Qin}
Let $G$ be a connected graph of order $n$. Then

\emph{(1)} $mc(G) \neq \binom{n}{2}-1$, $mc(G) \neq \binom{n}{2}-3$.

\emph{(2)} $mc(G)=\binom{n}{2}-2$ if and only if $G=K_n-e$.

\emph{(3)} $mc(G)=\binom{n}{2}-4$ if and only if $G\in \{K_n-2K_2, K_n-P_3\}$.
\end{cor}

\subsection{Results for graph products}

Mao, Wang, Yanling and Ye \cite{Mao} studied the monochromatic connection numbers of
the following graph products.
\begin{thm}\cite{Mao}
Let $G$ and $H$ be connected graphs.

$(1)$ If neither $G$ nor $H$ is a tree, then
$$
\max\{|E(G)||V(H)|,|E(H)||V(G)|\}+2\leq {\rm mc}(G\Box H)\leq
|E(G)||V(H)|+(|E(H)|-1)|V(G)|+1.
$$

$(2)$ If $G$ is not a tree and $H$ is a tree, then
$$
|E(H)||V(G)|+2\leq {\rm mc}(G\Box H)\leq|E(G)||V(H)|+1.
$$

$(3)$ If both $G$ and $H$ are trees, then
$$
|E(G)||E(H)|+1\leq {\rm mc}(G\Box H)\leq |E(G)||E(H)|+2.
$$

Moreover, the lower bounds are sharp.
\end{thm}

\begin{cor}\cite{Mao}
Let $G$ and $H$ be a connected graph.

$(1)$ If neither $G$ nor $H$ is a tree, then ${\rm mc}(G\Box H)\geq
\max\{{\rm mc}(G)|V(H)|+2,{\rm mc}(H)|V(G)|+2\}$.

$(2)$ If $G$ is not a tree and $H$ is a tree, then ${\rm mc}(G\Box
H)\geq {\rm mc}(H)|V(G)|+2$.

$(3)$ If both $G$ and $H$ are trees, then ${\rm mc}(G\Box H)\geq
{\rm mc}(G){\rm mc}(H)+1$.
\end{cor}

\begin{thm}\cite{Mao}
Let $G$ and $H$ be connected graphs, and let $G$ be noncomplete.

$(1)$ If neither $G$ nor $H$ is a tree, then
$$
|E(G)||V(H)|^2+2\leq {\rm mc}(G\circ H)\leq
|E(H)||V(G)|+|E(G)||V(H)|^2-|V(H)|+1.
$$

$(2)$ If $G$ is not a tree and $H$ is a tree, then
$$
|E(H)||V(G)|(|V(H)|+1)+2\leq {\rm mc}(G\circ
H)\leq|E(H)||V(G)|+|E(G)||V(H)|^2-|V(H)|+1.
$$

$(3)$ If $H$ is not a tree and $G$ is a tree, then
$$
|E(H)||V(G)|^2+2\leq {\rm mc}(G\circ
H)\leq|E(H)||V(G)|+|E(G)||V(H)|^2-|V(H)|+1.
$$

$(4)$ If both $G$ and $H$ are trees, then
$$
|E(H)||E(G)|(|V(H)|+1)+1\leq {\rm mc}(G\circ
H))\leq|E(H)||E(G)|(|V(H)|+1)+|V(H)|.
$$

Moreover, the lower bounds are sharp.
\end{thm}

\begin{cor}\cite{Mao}
Let $G$ and $H$ be connected graphs.

$(1)$ If neither $G$ nor $H$ is a tree, then ${\rm mc}(G\circ H)\geq
{\rm mc}(G)|V(H)|^2+2$.

$(2)$ If $G$ is not a tree and $H$ is a tree, then ${\rm mc}(G\circ
H)\geq {\rm mc}(H)|V(G)|(|V(H)|+1)+2$.

$(3)$ If $H$ is not a tree and $G$ is a tree, then ${\rm mc}(G\circ
H)\geq {\rm mc}(H)|V(G)|^2+2$.

$(4)$ If both $G$ and $H$ are trees, then ${\rm mc}(G\circ H))\geq
{\rm mc}(G){\rm mc}(H)(|V(H)|+1)+1$.
\end{cor}

\begin{thm}\cite{Mao}
Let $G$ and $H$ be a connected graph, and at least one of $G$ and $H$ is not a complete graph.

$(1)$ If neither $G$ nor $H$ is a tree, then
$$
{\rm mc}(G\boxtimes H)\geq
\max\{|E(G)||V(H)|+2|E(H)||E(G)|+2,|E(H)||V(G)|+2|E(H)||E(G)|+2\}
$$
and
$$
{\rm mc}(G\boxtimes H)\leq
|E(G)||V(H)|+|E(H)||V(G)|+2|E(H)||E(G)|-\min\{|V(H)|,|V(G)|\}+1.
$$

$(2)$ If $G$ not a tree and $H$ is a tree, then
$$
|E(H)||V(G)|+2|E(H)||E(G)|+2\leq {\rm mc}(G\boxtimes H)\leq
|E(G)||V(H)|+2|E(H)||E(G)|+1.
$$

$(3)$ If both $G$ and $H$ are trees, then
$$
3|E(H)||E(G)|+1\leq {\rm mc}(G\boxtimes H)\leq
3|E(H)||E(G)|+\min\{|V(G)|,|V(H)|\}.
$$
Moreover, the lower bounds are sharp.
\end{thm}
\begin{cor}\cite{Mao}
Let $G$ and $H$ be a connected graph.

$(1)$ If neither $G$ nor $H$ is a tree, then
$$
{\rm mc}(G\boxtimes H)\geq \max\{|{\rm mc}(G)||V(H)|+2|{\rm
mc}(H)||mc(G)|+2,|{\rm mc}(H)||V(G)|+2|{\rm mc}(H)||{\rm
mc}(G)|+2\}.
$$

$(2)$ If $G$ not a tree and $H$ is a tree, then
$$
{\rm mc}(G\boxtimes H)\geq |{\rm mc}(H)||V(G)|+2|{\rm mc}(H)||{\rm
mc}(G)|+2.
$$

$(3)$ If both $G$ and $H$ are trees, then
$$
{\rm mc}(G\boxtimes H)\geq 3|{\rm mc}(H)||{\rm mc}(G)|+1.
$$
\end{cor}

\begin{thm}\cite{Mao}
Let $G$ and $H$ be nonbipartite graphs. Then
$$
|E(H)||E(G)|+2\leq {\rm mc}(G\times H)\leq 2|E(H)||E(G)|+1.
$$

Moreover, the lower bounds are sharp.
\end{thm}
\begin{cor}\cite{Mao}
Let one of $G$ and $H$ be a non-bipartite connected graph. Then
$$
{\rm mc}(G\times H)\geq |{\rm mc}(H)||{\rm mc}(G)|+2.
$$
\end{cor}

As an application of the above results,
they also studied the following graph classes.

We call $P_n\Box P_m$ a \emph{two-dimensional grid graph},
where $P_n$ and $P_m$ are paths on $n$ and $m$ vertices, respectively.

\begin{pro}\cite{Mao}

$(1)$ For the network $P_n\Box P_m \ (n\geq 3, m\geq 2)$,
$$
{\rm mc}(P_n\Box P_m)=nm-n-m+2.
$$

$(2)$ For network $P_n\circ P_m \ (n\geq 4, m\geq 3)$,
$$
{\rm mc}(P_n\circ P_m)=m^2n-m^2-n+2.
$$
\end{pro}

An \emph{$n$-dimensional mesh} is the Cartesian product of $n$
linear arrays.
Particularly, two-dimensional grid graph is a
$2$-dimensional mesh.
An \emph{$n$-dimensional hypercube} is an $n$-dimensional mesh,
in which all the $n$ linear arrays are
of size $2$.

\begin{pro}\cite{Mao}

$(1)$ For $n$-dimensional mesh $P_{L_1}\Box P_{L_2}\Box \cdots \Box
P_{L_n} \ (n\geq 4)$,
$$
{\rm mc}(P_{L_1}\Box P_{L_2}\Box \cdots \Box P_{L_n})\geq
(2\ell_1\ell_2-\ell_1-\ell_2)(\ell_3\ell_4\cdots \ell_n)+2.
$$

$(2)$ For network $P_{L_1}\circ P_{L_2}\circ
\cdots \circ P_{L_n}$,
$${\rm mc}(P_{L_1}\circ P_{L_2}\circ \cdots \circ P_{L_n})\geq (\ell_1\ell^2_2+\ell_1\ell_2-\ell_1-\ell^2_2)(\ell_3\ell_4\cdots \ell_n)^2+2.$$
\end{pro}

An \emph{$n$-dimensional torus} is the Cartesian product of $n$
cycles $R_1,R_2,\cdots,R_n$ of size at least three.

\begin{pro}\cite{Mao}

$(1)$ For network $R_1\Box R_2\Box \cdots \Box R_n$, $n\geq 4$
$$
{\rm mc}(R_1\Box R_2\Box \cdots \Box R_n)\geq r_1r_2\cdots r_n+2.
$$
where $r_i$ is the order of $R_i$ and $3\leq i\leq
n$.

$(2)$ For network $R_1\circ R_2\circ \cdots \circ R_n$, $n\geq 4$
$$
{\rm mc}(R_1\circ R_2\circ \cdots \circ R_n)\geq r_1(r_2\cdots
r_n)^2+2.
$$
\end{pro}

Let $K_m$ be a clique of $m$ vertices, $m\geq 2$. An
\emph{$n$-dimensional generalized hypercube} is the Cartesian product of $n$ cliques.

\begin{pro}\cite{Mao}

$(1)$ For network $K_{m_1}\Box K_{m_2}\Box \cdots \Box K_{m_n} \
(m_i\geq 2, \ n\geq 3, \ 1\leq i\leq n)$
$$
{\rm mc}(K_{m_1}\Box K_{m_2}\Box \cdots \Box K_{m_n})\geq
{{m_1}\choose{2}}m_2\cdots m_n+2.
$$

$(2)$ For network $K_{m_1}\circ K_{m_2}\circ \cdots \circ K_{m_n}$,
$$
{\rm mc}(K_{m_1}\circ K_{m_2}\circ \cdots \circ
K_{m_n})={{m_1m_2\cdots m_n}\choose{2}}.
$$
\end{pro}

An \emph{$n$-dimensional hyper Petersen network}
$HP_n$ is the Cartesian product of $Q_{n-3}$ and
the well-known Petersen graph , where $n\geq 3$ and
$Q_{n-3}$ denotes an $(n-3)$-dimensional hypercube.

The network $HL_n$ is the lexicographical product of
$Q_{n-3}$ and the Petersen graph, where $n\geq 3$ and $Q_{n-3}$
denotes an $(n-3)$-dimensional hypercube.
\begin{pro}\cite{Mao}

$(1)$ For network $HP_3$ and $HL_3$, ${\rm mc}(HP_3)={\rm
mc}(HL_3)=7$;

$(2)$ For network $HL_4$ and $HP_4$, ${\rm mc}(HP_4)=22$ and
$112\leq {\rm mc}(HL_4)\leq 124$.
\end{pro}

For the join of two graphs, Jin, Li and Wang got the following results.
\begin{thm}\cite{Jin}
Let $G_1$ and $G_2$ be two disjoint connected graphs and $G =
G_1 + G_2$. Then $mc(G) = mc(G_1) + mc(G_2) + |V (G_1)||V (G_2)|$.
\end{thm}

\begin{thm}\cite{Jin}
Let $G$ be the join of a connected graph $G_1$ and a disconnected
graph $G_2$, where $V (G_1) \cap V (G_1) = \emptyset$. Then
$mc(G) = mc(G_1) + |E(G_2)| + |V (G_1)||V (G_2)| - |V (G_2)| + 1$.
\end{thm}

\begin{thm}\cite{Jin}
Let $G$ be the join of two disjoint disconnected graphs $G_1$
and $G_2$. Then $mc(G) = |E(G_1)| + |E(G_2)| + |V (G_1)||V (G_2)| - |V (G_1)| -
|V (G_2)| + 2$.
\end{thm}

\subsection{Results for random graphs}

The goal of $MC$-coloring of a graph is to find as many as colors to make the graph monochromatically connected.
So it is interesting to consider the threshold function of property $mc\left(G\left(n,p\right)\right)\ge f(n)$, where $f(n)$ is a function of $n$.
For any graph $G$ with $n$ vertices and any function $f(n)$, having $mc(G) \ge f(n)$ is a monotone graph property (adding
edges does not destroy this property), so it has a sharp threshold function.

Gu, Li, Qin and Zhao\cite{Qin} showed a sharp threshold function for $mc(G)$ as follows.

\begin{thm}\cite{Qin}
Let $f(n)$ be a function satisfying $1\leq f(n)<\frac{1}{2}n(n-1)$. Then
\begin{equation*}
p=
\left\{
  \begin{array}{ll}
   \frac{f(n)+n\log\log n}{n^2} & \hbox{ if $\ell n \log n\leq f(n)<\frac{1}{2}n(n-1)$, where $\ell\in \mathbb{R}^+$,} \\
   \frac{\log n}{n} & \hbox{ if  $f(n)=o(n\log n)$.} \end{array}
\right.
\end{equation*}
is a sharp threshold function for the property $mc\left(G\left(n,p\right)\right)\ge f(n)$.
\end{thm}
\noindent {\bf Remark 2.33.} Note that $mc\left(G\left(n,p\right)\right)\le \frac{1}{2}n(n-1)$ for any function $0\le p\le 1$, and $mc\left(G\left(n,p\right)\right)= \frac{1}{2}n(n-1)$ if and only if
$G(n, p)$ is isomorphic to the complete graph $K_n$. Hence we only concentrate on the case $f(n)<\frac{1}{2}n(n-1)$.

\section{The vertex-coloring version}

\subsection{Upper and lower bounds for $vmc(G)$}

For a connected graph $G$ of order 1 or 2, it is easy to check $vmc(G)=1, \ 2$, respectively.
For a connected graph $G$ of order at
least 3, Cai, Li and Wu \cite{Di2} got that a general lower bound for $vmc(G)$
is $\ell(T)+1\geq 3$, where $T$ is a spanning tree of $G$, and $\ell(T)$ is the number of leaves in $T$.
Simply take a spanning tree $T$ of $G$. Then, give all the
non-leaves in $T$ one color, and each leaf in $T$ a distinct new
color. Clearly, this is a VMC-coloring of $G$ using $\ell(T)+1$ colors.

By the known results about spanning trees with many leaves in
\cite{Caro2, Griggs, Kleitman}, Cai, Li and Wu \cite{Di2} got the following lower bounds.
\begin{pro}\cite{Di2}
Let $G$ be a connected graph with $n$ vertices and minimum degree $\delta$.

$(1)$ If $\delta\geq3$, then $vmc(G)\geq \frac{1}{4}n+3$.

$(2)$ If $\delta\geq4$, then $vmc(G)\geq \frac{2}{5}n+\frac{13}{5}$.

$(3)$ If $\delta\geq5$, then $vmc(G)\geq \frac{1}{2}n+3$.

$(4)$ If $\delta\geq3$, then $vmc(G)\geq
\left(1-\frac{\ln(\delta+1)}{\delta+1}(1+o_{\delta}(1))\right)n+1$.
\end{pro}
They also got an upper bound for $vmc(G)$.
\begin{pro}\cite{Di2}
Let $G$ be a connected graph with $n$ vertices and diameter $d$.

$(1)$ $vmc(G)=n$ if and only if $d\leq 2$;

$(2)$ If $d\geq 3$, then $vmc(G)\leq n-d+2$ , and the bound is sharp.
\end{pro}

\subsection{Erd\H{o}s-Gallai-type problems for $vmc(G)$}

Cai, Li and Wu \cite{Di2} studied two Erd\H{o}s-Gallai-type problems
for the graph parameter $vmc(G)$.

\noindent {\bf Problem A:} Given two positive integers $n$, $k$ with
$3\leq k\leq n$, compute the minimum integer $f_v(n,k)$ such that
for any graph $G$ of order $n$, if $|E(G)|\geq f_v(n,k)$ then
$vmc(G)\geq k$.

\noindent {\bf Problem B:} Given two positive integers $n$, $k$
with $3\leq k\leq n$, compute the maximum integer $g_v(n,k)$ such
that for any graph $G$ of order $n$, if $|E(G)|\leq g_v(n,k)$ then
$vmc(G)\leq k$.

Note that $g_v(n,n)={n \choose 2}$, and $g_v(n,k)$ does not exist for $3\leq k\leq n-1$.
This is because for a star $S_n$ on $n$ vertices, we have $vmc(S_n)=n$.
For this reason, Cai, Li and Wu \cite{Di2} just studied {\bf Problem A}.
They got the value of $f_v(n,k)$.
\begin{thm}\label{cor1}\cite{Di2}
Given two integers $n,k$ with $3\leq k\leq n$,
\begin{align*}
f_v(n,k)=
\begin{cases}
n-1 & \ $if$ \ \ k=3\\
n+\binom{k-2}{2} & \ $if$ \ \ 4\leq k \leq n-2\\
n-1+\binom{k-2}{2} & \ $if$ \ \ n-1\leq k \leq n\\
\end{cases}
\end{align*}
\end{thm}

\subsection{Nordhaus-Gaddum-type theorem for $vmc(G)$}

Cai, Li and Wu \cite{Di2} got the following Nordhaus-Gaddum-type result for $vmc(G)$
\begin{thm}\cite{Di2}
Let $G$ be a connected graph on $n\geq 5$ vertices
with connected complement $\overline{G}$.
Then $n+3\leq vmc(G)+vmc(\overline{G})\leq 2n$,
and $3n\leq vmc(G)\cdot vmc(\overline{G})\leq n^2$.
Moreover, these bounds are sharp.
\end{thm}

\section{The arc-coloring version for digraphs}

Gonz¨¢lez-Moreno, Guevara, and Montellano-Ballesteros \cite{Gonz¨¢lez} got
the following result for strongly connected oriented graph.
\begin{thm}\label{D-thm1}\cite{Gonz¨¢lez}
Let $D$ be a strongly connected oriented graph of size $m$,
and let $\Omega(D)$ be the minimum size of a strongly connected spanning subdigraph of $D$.
Then
\begin{align*}
smc(D)=m-\Omega(D)+1.
\end{align*}
\end{thm}

As an application of Theorem \ref{D-thm1},
they found a sufficient and necessary condition to determine whether a  strongly connected oriented graph is Hamiltonian.
\begin{cor}\label{D-cor1}\cite{Gonz¨¢lez}
Let $D$ be a strongly connected oriented graph of size $m$ and order $n$.
Then $D$ is Hamiltonian if and only if $smc(D)=m-n+1$.
\end{cor}

From Corollary \ref{D-cor1}, one can see that computing $\Omega(D)$ is NP-hard.

\section{Monochromatic indices}

\subsection{Edge version}

Li and Wu \cite{Di1} completely determined the $k$-monochromatic index for $k\geq 3$.
\begin{thm}\cite{Di1}
Let $G$ be a connected graph with $n$ vertices and $m$ edges.
Then $mx_k(G)=m-n+2$ for each $k$ with $3\leq k\leq n$.
\end{thm}

\subsection{Vertex version}

Li and Wu \cite{Di1} studied the hardness for computing $vmx_k(G)$.
They showed that given a connected graph $G=(V,E)$, and a positive integer $L$ with $L\leq |V|$,
to decide whether $vmx_k(G)\geq L$ is NP-complete for each $k$ with $2\leq k\leq |V|$.
In particular, computing $vmx_k(G)$ is NP-hard.

\subsection{Nordhaus-Gaddum-type results}

Recall that Cai, Li and Wu \cite{Di2} got the Nordhaus-Gaddum-type result for $vmc(G)$.
Li and Wu \cite{Di1} got the following Nordhaus-Gaddum-type lower bounds of $vmx_k$
for $k$ with $3\leq k\leq n$.

\begin{thm}\cite{Di1}
Suppose that both $G$ and $\overline{G}$ are connected graphs on $n$ vertices.
For $n=5$, $vmx_k(G)+vmx_k(\overline{G})\geq 6$ for $k$ with $3\leq k\leq 5$.
For $n=6$, $vmx_k(G)+vmx_k(\overline{G})\geq 8$ for $k$ with $3\leq k\leq 6$.
For $n\geq 7$,
if $n$ is odd, then $vmx_k(G)+vmx_k(\overline{G})\geq n+3$ for $k$ with $3\leq k\leq \frac{n-1}{2}$,
and $vmx_k(G)+vmx_k(\overline{G})\geq n+2$ for $k$ with $\frac{n+1}{2}\leq k\leq n$;
if $n=4t$, then $vmx_k(G)+vmx_k(\overline{G})\geq n+3$ for $k$ with $3\leq k\leq \frac{n}{2}-1$,
and $vmx_k(G)+vmx_k(\overline{G})\geq n+2$ for $k$ with $\frac{n}{2}\leq k\leq n$;
if $n=4t+2$, then $vmx_k(G)+vmx_k(\overline{G})\geq n+3$ for $k$ with $3\leq k\leq \frac{n}{2}$,
and $vmx_k(G)+vmx_k(\overline{G})\geq n+2$ for $k$ with $\frac{n}{2}+1\leq k\leq n$.
Moreover, all the above bounds are sharp.
\end{thm}

They also got the following Nordhaus-Gaddum-type upper bound of $vmx_k$
for $k$ with $\lceil\frac{n}{2}\rceil\leq k\leq n$.
\begin{thm}\cite{Di1}
Suppose that both $G$ and $\overline{G}$ are connected graphs on $n\geq 5$ vertices.
Then, for any $k$ with $\lceil\frac{n}{2}\rceil\leq k\leq n$, we have that
$vmx_k(G)+vmx_k(\overline{G})\leq 2n-2$,
and this bound is sharp.
\end{thm}

\section{The total-coloring version}

Jiang, Li and Zhang \cite{Jiang2} studied the hardness for computing $tmc(G)$.
They showed that
given a connected graph $G=(V,E)$, and a positive integer $L$ with $L\leq |V|+|E|$,
to decide whether $tmc(G)\geq L$ is NP-complete.
In particular, computing $tmc(G)$ is NP-hard.

\subsection{Upper and low bounds for $tmc(G)$}

Let $l(T)$ denote the number of leaves in a tree $T$.
For a connected graph $G$, let $l(G)=\max\{ \ l(T) \ |  $ $T$ is a spanning tree of $G$ $\}$.
Jiang, Li and Zhang \cite{Jiang1} got the following lower bound of $tmc(G)$.
\begin{thm}\cite{Jiang1}
 For a connected graph $G$ of order $n$ and size $m$, we have
$\mathrm{tmc}(G)\geq m-n+2+l(G)$.
\end{thm}

They also gave some sufficient conditions for graphs attaining this lower bound.
\begin{thm}\cite{Jiang1}\label{Jiang-thm3}
Let $G$ be a connected graph of order $n>3$ and size $m$. If $G$ has any of the following properties,
then $\mathrm{tmc}(G)=m-n+2+l(G)$.

$(a)$ The complement $\overline{G}$ of $G$ is $4$-connected.

$(b)$ $G$ is $K_3$-free.

$(c)$ $\Delta(G)<n-\frac{2m-3(n-1)}{n-3}$.

$(d)$ $diam(G)\geq 3$.

$(e)$ $G$ has a cut vertex.
\end{thm}

The upper bound of $\Delta(G)$ in Theorem\ref{Jiang-thm3}(c) is best possible. For example, let $G=K_{n-2,1,1}$.
Then $\mathrm{tmc}(G)=m-n+3+l(G)$ and $\Delta(G)=n-1=n-\frac{2m-3(n-1)}{n-3}$.

Jiang, Li and Zhang \cite{Jiang1} computed the total monochromatic connection numbers of
wheel graphs and complete multipartite graphs.
\begin{pro}\cite{Jiang1}
Let $G$ be a wheel $W_{n-1}$ of order $n\geq 5$ and size $m$. Then $\mathrm{tmc}(G)=m-n+2+l(G)$.
\end{pro}
\begin{pro}\cite{Jiang1}
Let $G= K_{n_1,\ldots,n_r}$ be a complete multipartite graph with $n_1 \geq \ldots \geq n_t\geq 2$ and $n_{t+1}=\ldots=n_r=1$.
Then $\mathrm{tmc}(G)=m+r-t$.
\end{pro}

\subsection{Comparing $\mathrm{tmc}(G)$ with $\mathrm{vmc}(G)$ and $\mathrm{mc}(G)$ }

Jiang, Li and Zhang \cite{Jiang1} compared $\mathrm{tmc}(G)$ with $\mathrm{vmc}(G)$
from different aspects.
\begin{thm}\cite{Jiang1}\label{Jiang-thm1}
Let $G$ be a connected graph of order $n$, size $m$ and diameter $d$. If $m\geq 2n-d-2$, then $\mathrm{tmc}(G)>\mathrm{vmc}(G)$.
\end{thm}
\begin{thm}\cite{Jiang1}\label{Jiang-thm2}
Let $G$ be a connected graph of order $n$, diameter $2$ and maximum degree $\Delta$. If $\Delta\geq \frac{n+1}{2}$, then $\mathrm{tmc}(G)>\mathrm{vmc}(G)$.
\end{thm}
Note that $\mathrm{tmc}(C_5)=4<\mathrm{vmc}(C_5)=5$, where $m<2n-d-2$ and $\Delta < \frac{n+1}{2}$. This implies that the conditions of Theorems \ref{Jiang-thm1} and \ref{Jiang-thm2} cannot be improved.
If $G$ is a star, then $\mathrm{tmc}(G)=\mathrm{vmc}(G)=n$.
However, they could not show whether there exist other graphs with $\mathrm{tmc}(G)\leq \mathrm{vmc}(G)$.
Then they proposed the following problem.

\begin{prob}\cite{Jiang1}
 Dose there exists a graph of order $n\geq 6$ except for the star graph such that $\mathrm{tmc}(G)\leq \mathrm{vmc}(G)$?
\end{prob}

In addition, they proposed the following conjecture.
\begin{con}\cite{Jiang1}
For a connected graph $G$, it always holds that $\mathrm{tmc}(G)>\mathrm{mc}(G)$.
\end{con}

Finally, they compared $tmc(G)$ with $mc(G)+vmc(G)$.
\begin{thm}\cite{Jiang1}
Let $G$ be a connected graph. Then $\mathrm{tmc}(G)\leq \mathrm{mc}(G)+\mathrm{vmc}(G)$, and the equality holds if and only if $G$ is a complete graph.
\end{thm}

\subsection{Results for graph classes}

Jiang, Li and Zhang \cite{Jiang2} characterized all connected graphs $G$ of order $n$ and size $m$
with $tmc(G)\in\{3,4,5,6,m+n-2,m+n-3,m+n-4\}$, respectively.
Let $\mathcal{T}_i$ denote the set of the trees with $l(G)=i$, where $2\leq i\leq n-1$. Note that if $G$ is a connected graph with $l(G)=2$, then $G$ is either a path or a cycle.

\begin{figure}[htbp]
\begin{center}
\scalebox{0.5}[0.5]{\includegraphics{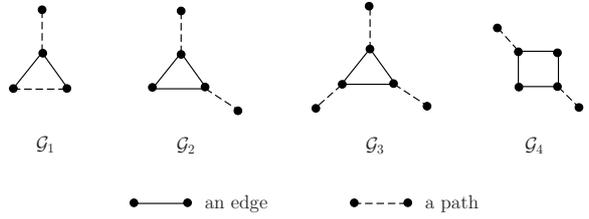}}
\end{center}
\caption{Unicyclic graphs with $l(G)=3$.}\label{Fig.2.}
\end{figure}

\begin{figure}[htbp]
\begin{center}
\scalebox{0.5}[0.5]{\includegraphics{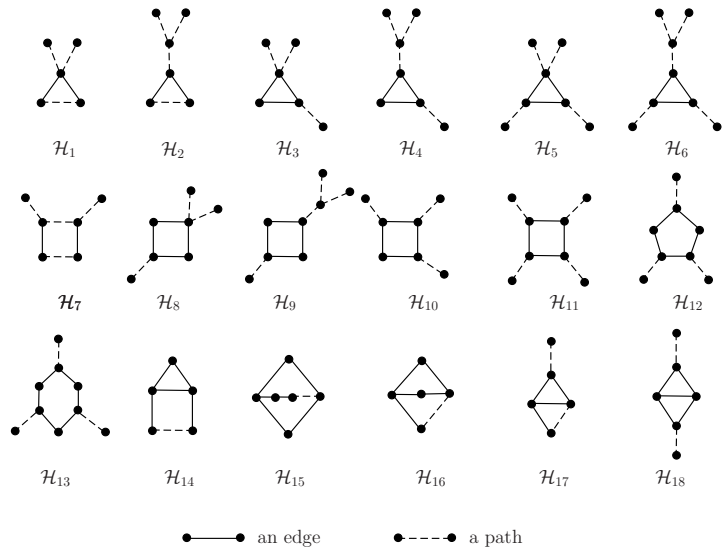}}
\end{center}
\caption{The graphs in Theorem \ref{Jiang-thm4}.}\label{Fig.3.}
\end{figure}

\begin{thm}\cite{Jiang2}
Let $G$ be a connected graph. Then $tmc(G)=3$ if and only if $G$ is a path.
\end{thm}
\begin{thm}\cite{Jiang2}
Let $G$ be a connected graph. Then $tmc(G)=4$ if and only if $G\in\mathcal{T}_3$ or $G$ is a cycle except for $K_3$.
\end{thm}
\begin{thm}\cite{Jiang2}
 Let $G$ be a connected graph. Then $tmc(G)=5$ if and only if $G\in\mathcal{T}_4$ or $G\in\mathcal{G}_i$, where $1\leq i\leq4$; see Figure \ref{Fig.2.}.
\end{thm}
\begin{thm}\label{Jiang-thm4}\cite{Jiang2}
Let $G$ be a connected graph. Then $tmc(G)=6$ if and only if $G=K_3$, $G\in\mathcal{T}_5$ or $G\in\mathcal{H}_i$, where $1\leq i\leq18$; see Figure \ref{Fig.3.}.
\end{thm}
\begin{thm}\cite{Jiang2}
Let $G$ be a connected graph. Then $tmc(G)=m+n-2$ if and only if $G=K_n-K_2$.
\end{thm}
\begin{thm}\cite{Jiang2}
 Let $G$ be a connected graph. Then $tmc(G)=m+n-3$ if and only if $G$ is either $K_n-K_3$ or $K_n-P_3$.
\end{thm}
\begin{thm}\cite{Jiang2}
Let $G$ be a connected graph. Then $tmc(G)=m+n-4$ if and only if $G\in\{K_n-P_4,K_n-2K_2,K_n-K_4,K_n-(K_4-K_2),K_n-(K_4-P_3),K_n-C_4,K_n-K_{1,3}\}$.
\end{thm}

\subsection{Results for random graphs}

For a property $P$ of graphs and a positive integer $n$, define $Prob(P,n)$ to be the ratio of the number of graphs with $n$
labeled vertices having property $P$ over the total number of graphs with these vertices. If $Prob(P,n)$ approaches 1 as $n$ tends
to infinity, then we say that {\it almost all} graphs have property $P$. More details can be found in \cite{BH}.
Jiang, Li and Zhang \cite{Jiang1} got the following result for $tmc(G)$.
\begin{thm}\cite{Jiang1} For almost all graphs $G$ of order $n$ and size $m$, we have $\mathrm{tmc}(G)=m-n+2+l(G)$.
\end{thm}
Jiang, Li and Zhang \cite{Jiang2} showed a sharp threshold function for $tmc(G)$ as follows.
\begin{thm}\cite{Jiang2}
Let $f(n)$ be a function satisfying $1\leq f(n)<\frac{1}{2}n(n-1)+n$. Then
\begin{eqnarray}p=
\begin{cases}
\frac{f(n)+n\log\log n}{n^2} &if\ ln\log n\leq f(n)< \frac{1}{2}n(n-1)+n, \cr
 \  & where\ l\in \mathbb{R}^{+}, \cr
\frac{\log n}{n} &if\ f(n)=o(n\log n).
\end{cases}
\end{eqnarray}
is a sharp threshold function for the property $tmc(G(n,p))\geq f(n)$.
\end{thm}

\noindent {\bf Remark 6.19.} Note that if $f(n)=\frac{1}{2}n(n-1)+n$, then $G(n,p)$ is a complete graph $K_n$ and $p=1$.
Hence we only concentrate on the case $f(n)<\frac{1}{2}n(n-1)+n$.

\subsection{Erd\H{o}s-Gallai-type problems for $tmc(G)$}

Jiang, Li and Zhang \cite{Jiang3} studied the following two kinds of Erd\H{o}s-Gallai-type problems for $tmc(G)$.

\noindent\textbf{Problem A.}
Given two positive integers $n$ and
$k$ with $3\leq k\leq {n \choose 2}+n$, compute the minimum integer $f_T(n,k)$ such that
for any graph $G$ of order $n$, if $|E(G)|\geq f_T(n,k)$ then
$tmc(G)\geq k$.

\noindent\textbf{Problem B.} Given two positive integers $n$ and $k$
with $3\leq k\leq {n \choose 2}+n$, compute the maximum integer $g_T(n,k)$ such
that for any graph $G$ of order $n$, if $|E(G)|\leq g_T(n,k)$ then
$tmc(G)\leq k$.

They completely determined the values of $f_T(n,k)$ and $g_T(n,k)$.

\begin{thm}\cite{Jiang3}
Given two positive integers $n$ and $k$ with $3\leq k\leq\binom{n}{2}+n$,
\begin{eqnarray}f_T(n,k)=
\begin{cases}
n-1 &if\ k=3, \cr
n+k-t-2 &if\ k=\binom{t}{2}+t+2-s,\ where\ 0\leq s\leq t-1 \ and\ 2\leq t\leq n-2, \cr
k &if\ \binom{n}{2}-n+4\leq k\leq \binom{n}{2}+n-3\lfloor\frac{n}{2}\rfloor\ except\ for\ n\
is\ odd\cr
\ \ &\ \ \ \ and\ k=\binom{n}{2}+n-3\lfloor\frac{n}{2}\rfloor,\cr
\binom{n}{2}-r &if\ \binom{n}{2}+n-3(r+1)<k\leq\binom{n}{2}+n-3r,\ where\ 0\leq r\leq \lfloor\frac{n}{2}\rfloor-1\cr
\ \ &\ \ \ \ or\ n\
is\ odd,\ r=\lfloor\frac{n}{2}\rfloor\ and\ k=\binom{n}{2}+n-3\lfloor\frac{n}{2}\rfloor.
 \end{cases}
\end{eqnarray}
\end{thm}

\begin{thm}\cite{Jiang3}
Given two positive integers $n$ and $k$ with $n\leq k\leq \binom{n}{2}+n$,
\begin{eqnarray}g_T(n,k)=
\begin{cases}
k-n+t &if\ \binom{n-t}{2}+t(n-t-1)+n\leq k\leq \binom{n-t}{2}+t(n-t)+n-2, \cr
k-n+t-1 &if\ k=\binom{n-t}{2}+t(n-t)+n-1, \cr
\binom{n}{2}-1 &if\ k=\binom{n}{2}+n-1, \cr
\binom{n}{2} &if\ k=\binom{n}{2}+n,
\end{cases}
\end{eqnarray}
for $2\leq t\leq n-1$.
\end{thm}

\section{Concluding remarks}

This survey tries to summarize all the results on monochromatic connection of graphs in the existing literature. The simple purpose is to promote
the research along this subject. As one can see, there are some basic problems remaining unsolved. For example,
what is the computational complexity of determining the monochromatic connection number $mc(G)$ for a given connected graph
$G$ ? From Theorem 2.1 (d) one can see that this problem is reduced to only considering those graphs with diameter 2.
It is easily seen also from Theorem 2.1 (a) that for almost all connected graphs $G$ it holds that $mc(G)=m(G)-n(G)+2$.

Another problem is to consider more monochromatic paths connecting a pair of vertices. The definitions can be easily
given as follows. An edge-colored graph is called monochromatically $k$-connected if each pair of vertices of the graph is connected by
$k$ monochromatic paths in the graph. For a $k$-connected graph $G$, the {\it monochromatic $k$-connection number}, denoted by $mc_k(G)$, is defined as the
maximum number of colors that are needed in order to make $G$ monochromatically $k$-connected. As far as we knew, there is no paper published
on this parameter. We think that to get some bounds for the case $k=2$ is already quite interesting and not so easy.

It is seen that results for the monochromatic indices are very few, and more efforts are needed for deepening the research. It is
also seen that research on $smc(D)$ for digraphs has just started, and one can develop it with many possibilities.

Finally, we point out that we changed some terminology and notation. For examples, we use $vmc(G)$ to replace $mvc(G)$ and $vmx_k(G)$ to replace $mvx_k(G)$,
etc. This is because we think that the term ``vertex-monochromatic connection" is better than ``monochromatic vertex-connection".
This is just a matter of taste, depending on authors and readers.


\begin{thebibliography}{111}

\bibitem{BH}
A. Blass, F. Harary, \emph{Properties of almost all graphs and complexes}, J. Graph Theory 3(3)(1979), 225-240.

\bibitem{BT}
B. Bollob\'{a}s, A. Thomason, \emph{Threshold functions}, Combinatorica 7(1986), 35--38.

\bibitem{Bondy} J.A. Bondy and U.S.R. Murty, \emph{Graph Theory with Applications}, The Macmillan Press, London and Basingstoker, 1976.

\bibitem{Chartrand1}  G. Chartrand, G. Johns, K. McKeon, P. Zhang,
\emph{Rainbow connection in graphs}, Math. Bohem. 133(2008), 85-98.

\bibitem{Cai} Q. Cai, X. Li, D. Wu, \emph{Erd\H{o}s-Gallai-type results for colorful
monochromatic connectivity of a graph}, J. Comb. Optim. 33(1)(2017), 123-131.

\bibitem{Di2} Q. Cai, X. Li, D. Wu, \emph{Some extremal results on the colorful monochromatic vertex-connectivity of a graph}, arXiv:1503.08941.

\bibitem{Caro1}  Y. Caro, A. Lev, Y. Roditty, Zs. Tuza, R. Yuster,
\emph{On rainbow connection}, Electron. J. Combin. 15(1)(2008), R57.

\bibitem{Caro2} Y. Caro, D.B. West, R. Yuster, \emph{Connected domination and spanning
trees with many leaves}, SIAM J. Discrete Math. 13(2)(2000), 202-211.

\bibitem{Caro} Y. Caro, R. Yuster, \emph{Colorful monochromatic connectivity}, Discrete Math. 311(2011), 1786-1792.

\bibitem{ER}
P. Erd\"{o}s, A. R\'{e}nyi, \emph{On the evolution of random graphs}, Publ. Math. Inst. Hungar. Acad. Sci. 5(1960), 17--61.

\bibitem{FK}
E. Friedgut, G. Kalai, \emph{Every monotone graph property has a sharp threshold},
Proc. Amer. Math. Soc. 124(1996), 2993--3002.

\bibitem{Gonz¨¢lez}
D. Gonz¨¢lez-Moreno, M. Guevara, J.J. Montellano-Ballesteros, \emph{Monochromatic connecting colorings in strongly connected oriented graphs},
Discrete Math. 340(4)(2017), 578-584.

\bibitem{Griggs} J.R. Griggs, M. Wu, \emph{Spanning trees in graphs of
minimum degree 4 or 5}, Discrete Math. 104(2)(1992), 167-183.

\bibitem{Qin} R. Gu, X. Li, Z. Qin, Y. Zhao, \emph{More on the colorful monochromatic connectivity},
Bull. Malays. Math. Sci. Soc. DOI 10.1007/s40840-015-0274-2, in press.

\bibitem{Jiang1}
H. Jiang, X. Li, Y. Zhang, \emph{Total monochromatic connection of graphs}, Discrete Math. 340(2017), 175-180.

\bibitem{Jiang2}
H. Jiang, X. Li, Y. Zhang, \emph{More on total monochromatic connection of graphs}, Ars Combin. 136(2018).

\bibitem{Jiang3}
H. Jiang, X. Li, Y. Zhang, \emph{Erd\H{o}s-Gallai-type results for total monochromatic connection of graphs}, arXiv:1612.05381.

\bibitem{Jin} Z. Jin, X. Li, K. Wang, The monochromatic connectivity of some graphs, Manuscript, 2016.

\bibitem{Kleitman} D.J. Kleitman, D.B. West,
\emph{Spanning trees with many leaves,} SIAM J. Discrete Math. 4(1)(1991), 99-106.

\bibitem{Krivelevich and Yuster} M. Krivelevich, R. Yuster,
\emph{The rainbow connection of a graph is (at most) reciprocal to its minimum degree},
J. Graph Theory. 63(3)(2010), 185-191.

\bibitem{Sun}  X. Li, Y. Shi, Y. Sun, \emph{Rainbow connections of graphs:
A survey}, Graphs \& Combin. 29(2013), 1-38.

\bibitem{LiSun2} X. Li, Y. Sun, \emph{Rainbow Connections of Graphs},
 SpringerBriefs in Math. Springer, New York, 2012.

\bibitem{Di1} X. Li, Di Wu, \emph{The (vertex-)monochromatic index of a graph},
J. Comb. Optim. 33(2017), 1443-1453.

\bibitem{Mao}
Y. Mao, Z. Wang, F. Yanling, C. Ye, \emph{Monochromatic connectivity and graph products},
Discrete Math, Algorithm. Appl. 8(01)2016, 1650011.


\end{thebibliography}
\end{document}